\documentclass[11pt]{article}
\input epsf
\usepackage{amsfonts}
\usepackage{amscd}
\usepackage{amsmath,amssymb}
\usepackage{amstext}
\usepackage{amscd}
\usepackage[all]{xy}
\usepackage{pstricks,pst-node,pst-tree}

\newtheorem{lem}{LEMMA}[section]
\newtheorem{theo}[lem]{THEOREM}

\newtheorem{coro}[lem]{COROLLARY}

\newtheorem{definition}[lem]{Definition}

\newtheorem{rem}[lem]{Remark}
\newtheorem{ex}[lem]{Example}
\newtheorem{exs}[lem]{Examples}
\renewcommand{\descriptionlabel}[1]%
       {\hspace{\labelsep}\textsf{#1}}

\begin{document}

\title{ Homogeneity of dynamically defined wild knots}
\author{Gabriela Hinojosa\thanks{This work was partially supported by PROMEP (SEP) and CONACyT (Mexico), 
grant G36357-E}, Alberto Verjovsky\thanks{This work was partially supported by
CONACyT (Mexico), grant G36357-E,
a joint project CNRS-CONACyT and   PAPIIT (UNAM) grant ES11140}}
\date{April 22, 2005}

\maketitle

\begin{abstract} In this paper we prove that a wild knot $K$ which is the
limit set of a Kleinian group acting conformally on the unit 3-sphere,
with its standard metric, is homogeneous: given two points
$p,\,\,q\in{K}$ there exists a homeomorphism $f$ of the sphere
such that $f(K)=K$ and $f(p)=q$. We also show that if the wild knot
is a fibered knot then we can choose an $f$ which preserves the fibers.

\end{abstract}

\section{Introduction}

The birth of wild topology was in the 1920's with works of Alexander,
Antoine, Artin and Fox, among others. At that time one of the main problems
was to generalize the Schoenflies Theorem. Let $S$ be a simple closed
surface in $\mathbb{R}^{3}$ which is homeomorphic to the unit sphere
$\mathbb{S}^{2}$. Let $h$ be a homeomorphism of $S$ onto the unit
sphere $\mathbb{S}^{2}$ in $\mathbb{R}^{3}$. Is there an extension
$\tilde{h}$ of $h$ such that $\tilde{h}$ is a homeomorphism of
$\mathbb{R}^{3}$ onto itself? Alexander proved this result in the special case
that $S$ is a finite polytope. At the same time, however, he gave his
famous example, the {\it Alexander horned sphere}, where its unbounded
complement in $\mathbb{R}^{3}$ is not simply connected and, in fact, its fundamental
group is infinitely generated. Since the
complement of $\mathbb{S}^{2}$ in $\mathbb{R}^{3}$ is not simply
connected, it follows that no homeomorphism of $\mathbb{R}^{3}$ onto
itself will send the horned sphere onto $\mathbb{S}^{2}$ (see
\cite{rolfsen}, \cite{hocking}). The work of Alexander was published
in 1924.

In 1948, Artin and Fox gave the definition of {\it tame embeddings}
and {\it wild embeddings} and constructed a number of surprising
examples. For instance, see Figure 1.\\

\begin{figure}[tbh]
\centerline{\epsfxsize=2in \epsfbox{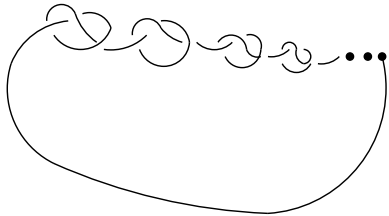}}
\caption{\sl An example of a wild knot.}
\label{F1}
\end{figure}

Many works by Antoine, Bing, Harold, Moise, Mazur, Brown, Montesinos, among
others have contributed much to the understanding of wild sets in
$\mathbb{R}^{3}$.

Let $K\subset\mathbb{S}^{3}$ be a knot. We say that a point $x\in K$ is 
{\it locally flat} if
 there exists an open neighborhood $U$ of $x$ such that
there is a homeomorphism of pairs: 
$(U,U\cap K)\sim (Int(\mathbb{B}^{3}),Int(\mathbb{B}^{1}))$.
Otherwise, $x$ is said to be a {\it wild} point.
A knot $K$ is a wild knot if it contains at least one wild point.

We say that a knot $K\subset\mathbb{S}^{3}$ is {\it homogeneous} 
if given two points
$p$, $q\in K$, there exists a homeomorphism 
$\psi:\mathbb{S}^{3}\rightarrow\mathbb{S}^{3}$
such that $\psi(K)=K$ and $\psi(p)=q$.

The wild knot $K$ given by Artin and Fox (Figure 1) is not
homogeneous. In fact, it contains just one wild point $p$, 
hence it is not possible to give a homeomorphism 
$\psi:\mathbb{S}^{3}\rightarrow\mathbb{S}^{3}$
such that $\psi(K)=K$ and  $\psi(p)=q$, $q\neq{p}$, 
since any homeomorphism sends
wild points into wild points. In general, wild knots are 
not homogeneous.

The purpose of this paper is to show that dynamically defined 
wild knots (see section 2) are homogeneous. In section 3, 
we will give a proof of this fact.
   
\section{Preliminaries}

In this section, we will describe the construction of dynamically defined wild knots.
We will begin with some basic definitions.

Let $M\ddot{o}b(\mathbb{S}^{n})$ denote the group of M\"obius
transformations of the n-sphere 
$\mathbb{S}^{3}=\mathbb{R}^{n}\cup\{\infty\}$ i.e. the group of diffeomorphisms of the
n-sphere that preserves angles with respect to the standard metric. Let 
$\Gamma\subset M\ddot{o}b(\mathbb{S}^{n})$ be a discrete subgroup. Then 
$x\in\mathbb{S}^{3}$ is a point of discontinuity for $\Gamma$ if there is a neighborhood $U$ of $x$ such
that $U\cap gU\neq\emptyset$ only for finitely many $g\in\Gamma$. The 
{\it domain of discontinuity} $\Omega(\Gamma)$ consists of all points of discontinuity.
\begin{definition}
(\cite{kap1}) A Kleinian group is a subgroup of $M\ddot{o}b(\mathbb{S}^{n})$ with non-empty domain 
of discontinuity. The complement 
$\mathbb{S}^{n}-\Omega(\Gamma)=\Lambda(\Gamma)$ is called {\it  limit set} of $\Gamma$.
\end{definition}
Next, we will give the construction of dynamically defined wild knots.

\begin{definition}
A  necklace $T_{1}$ of $n$-pearls ($n\geq 3$), is a collection of $n$ consecutive 
2-spheres $\Sigma_{1}$, $\Sigma_{2},\ldots ,\Sigma_{n}$ in
$\mathbb{S}^{3}$, such that $\Sigma_{i}\cap\Sigma_{j}=\emptyset$
($j\neq i+1, i-1$ mod $n$), except that
$\Sigma_{i}$ and $\Sigma_{i+1}$ are tangent ($i=1,2,\ldots,n-1$) and
$\Sigma_{1}$ and $\Sigma_{n}$ are 
tangent. Each
2-sphere is called a {\it  pearl}. 
\end{definition}

If the points of tangency are
joined by spherical geodesic segments in $\mathbb{S}^{3}$, we obtain a 
polygonal knot $K_{1}$. It is called the {\it polygonal template} of $T_{1}$. 
We define the {\it filled-in $T$} as $|T_{1}|=\cup^{n}_{i=1} B_{i}$, 
where $B_{i}$ is the round closed 3-ball whose boundary $\partial
B_{i}$ is $\Sigma_{i}$.

\begin{ex}
 $K=$ Trefoil knot.
\begin{figure}[tbh]
\centerline{\epsfxsize=1in \epsfbox{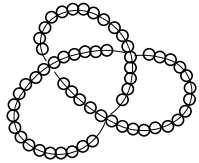}}
\caption{\sl A pearl-necklace whose template is the trefoil knot.}
\label{F2}
\end{figure}
\end{ex}

Let $\Gamma$ be the group generated by reflections $I_{j}$, through 
$\Sigma_{j}$ ($j=1,\ldots,n$). Then $\Gamma$ is a Kleinian group. 
We will describe geometrically the action of $\Gamma$.

\begin{enumerate}

\item First stage: Observe that when we reflect with respect to each $\Sigma_{k}$ 
($k=1,2,\ldots,n$), a mirror image of $K_{1}$ is mapped into the ball $B_{k}$ (see Figure 3). 

\begin{figure}[tbh]
\centerline{\epsfxsize=1in \epsfbox{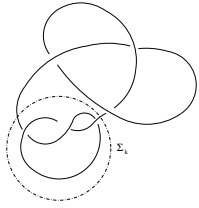}}
\caption{\sl Reflection with respect to $\Sigma_{k}$.}
\label{F3}
\end{figure}

After reflecting with respect to each pearl, we obtain a new necklace
$T_{2}$ of $n(n-1)$ pearls, subordinate to a new knot $K_{2}$; which
is in turn isotopic to the connected sum of $n+1$ copies of $K_{1}$ (see
Figure 2).

\begin{figure}[tbh]
\centerline{\epsfxsize=1in \epsfbox{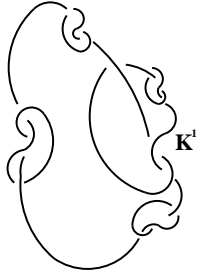}}
\caption{\sl A schematic figure of the reflecting process first step .}
\label{F4}
\end{figure}

\item Second stage: Now, we reflect with respect to each pearl of $T_{2}$. When
we are finished, we obtain  a new necklace $T_{3}$ of $n(n-1)^{2}$ pearls. Its
template is a polygonal knot $K_{3}$; which is in turn isotopic to the connected
sum of $n^{2}-n+1$ copies of $K_{1}$ and $n$ copies of its mirror image (recall that composition
of an even number of reflections is orientation-preserving). Observe that  $|T_{3}|\subset
|T_{2}|$ (see Figure 5). 

\begin{figure}[tbh]
\centerline{\epsfxsize=1in \epsfbox{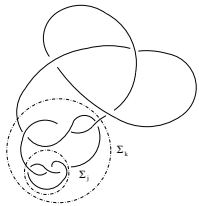}}
\caption{\sl Reflection with respect to $\Sigma_{k}$
         after reflecting with respect to $\Sigma_{j}$.}
\label{F5}
\end{figure}

\item $k^{th}$-Stage: We reflect with respect of each pearl of $T_{k}$.
At the end of this stage, we obtain a new  necklace  $T_{k+1}$ of $n(n-1)^{k}$ pearls, 
subordinate to a polygonal 
knot $K_{k+1}$. By construction, $|T_{k+1}|\subset |T_{k}|$.
\end{enumerate}

Then, the limit set is given by the inverse limit (see \cite{kap1}, \cite{maskit})
$$
\Lambda(\Gamma)=\varprojlim_{k} |T_{k}|=\bigcap_{k=1}^{\infty} |T_{k}|.
$$
It has been proved (see \cite{kap1}, \cite{maskit}) that the limit set $\Lambda(\Gamma)$ is 
a wild knot in the sense of Artin and Fox. It is called a {\it dynamically-defined wild knot.}

\section{Homogeneity}

Let $T_{1}$ be a $n_{1}$-pearl necklace subordinate to the polygonal knot $K_{1}$.
We can assume without loss of generality that 
$K_{1}\subset{\mathbb R}^3\subset{\mathbb S}^3=\mathbb{R}^{3}\cup\{\infty\}$.

Let $V_{1}$ be a closed
tubular neighborhood of  $K_{1}$ and $\pi_1:V_{1}\to{K_{1}}$ the projection.
We can assume that $\pi_1^{-1}(\{x\})$ is an euclidean 2-disk of radius $r_{1}>0$
independent of $x$. If $\{p_{1_{1}},\ldots, p_{1_{n_{1}}} \}$ are the points
of tangency of consecutive pearls of $T_{1}$, we can also assume 
that $\pi_1^{-1}(\{p_{1_{j}}\})$ is tangent to the consecutive pearls
at $p_{1_{j}}\,$ ($1\leq{j}\leq{n_{1}}$). The tubular neighborhood $V_{1}$
is the union of $n_{1}$ ``solid cylinders''  $V_{1_{1}}\ldots V_{1_{n_{1}}}$
where $V_{1_{j}}=\pi_1^{-1}(\{[p_{1_{j}},p_{1_{j+1}}]\})$ is called solid cylinder, since 
it is homeomorphic to a solid cylinder $C=\mathbb D^2\times[0,1]$ 
(see Figure 6).  

\begin{figure}[tbh]
\centerline{\epsfxsize=1.5in \epsfbox{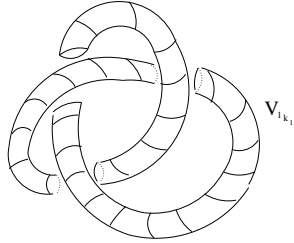}}
\caption{\sl A tubular neighbourhood as a union of ``cylinders''.}
\label{F6}
\end{figure}

For the second stage, we have a pearl necklace $T_{2}$ with $n_{2}$ pearls 
subordinate 
to the polygonal knot
$K_{2}$. Let $V_{2}\subset \mbox{Int}(V_{1})$ be a closed
tubular neighborhood of  $K_{2}$ and $\pi_2:V_{2}\to{K_{2}}$ be the projection.
We again assume that $\pi_2^{-1}(\{x\})$ is an euclidean 2-disk of radius $r_{2}>0$
independent of $x$. 
Notice that the points $\{p_{1_{1}},\ldots, p_{1_{n_{1}}} \}$ are also points
of tangency of consecutive pearls of  $T_{2}$. We will denote by 
$\{p_{1_{i},2_{1}},\ldots, p_{1_{i},2_{n-1}} \}\subset T_{2}$  the corresponding points
of tangency of consecutive pearls of $V_{1_{i}}\cap V_{2}$, $1\leq i \leq n$.
We can again assume 
that $\pi_2^{-1}(\{p_{1_{i},2_{j}}\})$ is tangent to the consecutive pearls
at $p_{1_{i},2_{j}}\,$ ($1\leq i \leq n$, $1\leq j\leq n-1$). The tubular neighborhood $V_{2}$
is the union of $n_{2}$ solid cylinders $V_{1_{i}2_{j}}$ ($1\leq i \leq n$, $1\leq j\leq n-1$)
where $V_{1_{i}2_{j}}=\pi_2^{-1}(\{[p_{1_{i}2_{j}},p_{1_{i}2_{j+1}}]\})$. 

We continue inductively, so
at the end of the $k$-stage of the 
reflecting process, we have the pearl necklace $T_{k}$ with $n_{k}$ pearls subordinate 
to the polygonal knot $K_{k}$.
Let $V_{k}$ be a closed
tubular neighborhood of  $K_{k}$ such that $V_{k}\subset \mbox{Int}(V_{k-1})$. Let
$\pi_k:V_{k}\to{K_{k}}$ be the projection.
We assume that $\pi_k^{-1}(\{x\})$ is an euclidean 2-disk of radius $r_{k}>0$
independent of $x$. We will denote by 
$\{p_{1_{i_{1}},2_{i_{2}},\ldots ,k_{1}},\ldots,p_{1_{i_{1}},2_{i_{2}},\ldots ,k_{n-1}} \}\subset T_{k}$  
the corresponding points
of tangency of consecutive pearls of 
$(V_{1_{i_{1}},2_{i_{2}},\ldots ,(k-1)_{i_{k-1}}})\cap V_{k}$, $1\leq i_{1} \leq n$
and $1\leq i_{2},\ldots i_{k-1} \leq n-1$. The tubular neighborhood $V_{k}$
is the union of $n_{k}$ solid cylinders $V_{1_{i_{1}},2_{i_{2}},\ldots ,(k-1)_{i_{k-1}},k_{i_{k}}}$.  
Notice that $\lim_{k\rightarrow\infty}r_{k}=0$ and 
$\Lambda=\cap_{k=1}^{\infty}V_{k}$.

Let $p$, $q\in\Lambda$. There exist two sequences of solid cylinders  
$\{V_{1_{i_{1}},2_{i_{2}},\ldots,n_{i_{n}}}\}$,  and  
$\{V_{1_{j_{1}},2_{j_{2}},\ldots,n_{j_{n}}}\}$ where $V_{1_{i_{1}},2_{i_{2}},\ldots,n_{i_{n}}}$ 
and $V_{1_{j_{1}},2_{j_{2}},\ldots,n_{j_{n}}}\in V_{n}$, such that 
$p=\cap_{n=1}^{\infty} V_{1_{i_{1}},2_{i_{2}},\ldots,n_{i_{n}}}$ and
  $q=\cap_{n=1}^{\infty} V_{1_{j_{1}},2_{j_{2}},\ldots,n_{j_{n}}}$. In fact,
these sequences converge to $p$ and $q$ 
respectively, with respect to the Hausdorff metric of closed sets on $\mathbb{S}^{3}$.

We define the homeomorphism 
$F_{0}:(\mathbb{S}^{3}-\mbox{Int}(V_{1}))\rightarrow (\mathbb{S}^{3}-\mbox{Int}(V_{1}))$ 
such that it sends $\partial V_{1_{i_{1}}}\cap\partial V_{1}$ into 
$\partial V_{1_{j_{1}}}\cap\partial V_{1}$, $\partial V_{1_{i_{1}+1}}\cap\partial V_{1}$ into
$\partial V_{1_{j_{1}+1}}\cap\partial V_{1}$ and so on. This map also sends
$\pi_{1}^{-1}(p_{1_{i_{1}}})\cap V_{1}$ into
$\pi_{1}^{-1}(p_{1_{j_{1}}})\cap V_{1}$, $\pi_{1}^{-1}(p_{1_{i_{1}+1}})\cap V_{1}$ into
$\pi_{1}^{-1}(p_{1_{j_{1}+1}})\cap V_{1}$ and so on. 

Next, we will define a homeomorphism 
$f_{1}:(\partial V_{1}\cup\partial V_{2})\rightarrow (\partial V_{1}\cup\partial V_{2})$ such that 
$f_{1}|_{\partial V_{1}}=F_{0}$.
Let $B_{1_{k_{1}}}=\cup_{k_{2}=1}^{n-1}V_{1_{k_{1}}2_{k_{2}}}$.
The pair  $(V_{1_{k_{1}}},B_{1_{k_{1}}}\cap V_{2})$
is a solid tangle; i.e. a solid cylinder with a knotted hole (see Figure 7), and 
it is homeomorphic to the solid tangle $(C,\tilde{K})$, where $C$ is a solid cylinder and 
$\tilde{K}$ is the mirror image of
the knot $K$ via the homeomorphism $h_{1_{k_{1}}}$, which will be used below. 

\begin{figure}[tbh]
\centerline{\epsfxsize=1.4in \epsfbox{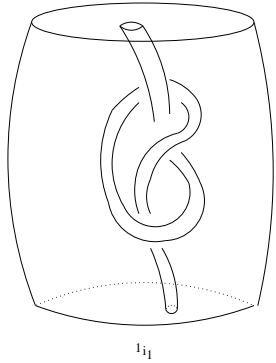}}
\caption{\sl A solid tangle.}
\label{F7}
\end{figure}

Since,
$F_{0}$ sends $\partial V_{1_{k_{1}}}\cap V_{1}$ into $\partial V_{1_{l_{1}}}\cap V_{1}$, 
we define the map $f_{1}$ in the following way. If
$k_{1}\neq i_{1}$, then $f_{1}$ sends the pair
$(\partial V_{1_{k_{1}}}, \partial B_{1_{k_{1}}}\cap V_{2})$ into 
$(\partial V_{1_{l_{1}}}, \partial B_{1_{l_{1}}}\cap V_{2})$, 
where $\partial V_{1_{k_{1}}2_{k_{2}}}\cap V_{2}$ is sent to $\partial V_{1_{l_{1}}2_{k_{2}}}\cap V_{2}$ and 
($\pi_{1}^{-1}(p_{1_{k_{1}}})-\mbox{Int}(\pi_{2}^{-1}(p_{1_{k_{1}}2_{k_{2}}}))$) is sent to
($\pi_{1}^{-1}(p_{1_{l_{1}}})-\mbox{Int}(\pi_{2}^{-1}(p_{1_{l_{1}}2_{k_{2}}}))$). If 
$k_{1}=i_{1}$, then $f_{1}$ sends
the pair
$(\partial V_{1_{i_{1}}}, \partial B_{1_{i_{1}}}\cap V_{2})$ into 
$(\partial V_{1_{j_{1}}}, \partial B_{1_{j_{1}}}\cap V_{2})$ such that
$\partial V_{1_{i_{1}}2_{i_{2}}}\cap V_{2}$ goes into $\partial V_{1_{j_{1}}2_{j_{2}}}\cap V_{2}$ and
($\pi_{1}^{-1}(p_{1_{i_{1}}})-\mbox{Int}(\pi_{2}^{-1}(p_{1_{i_{1}}2_{i_{2}}}))$) is sent to
($\pi_{1}^{-1}(p_{1_{j_{1}}})-\mbox{Int}(\pi_{2}^{-1}(p_{1_{j_{1}}2_{j_{2}}}))$) (see Figure 8).

\begin{figure}[tbh]
\centerline{\epsfxsize=3in \epsfbox{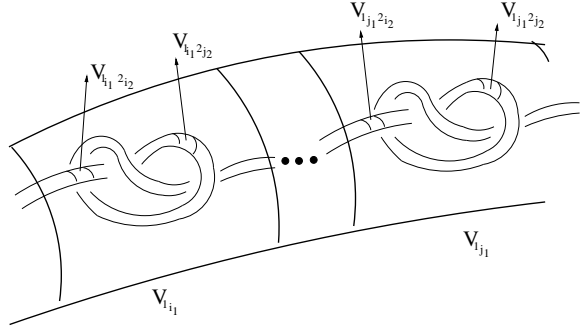}}
\caption{\sl Geometric description of the map $f_{1}$.}
\label{F8}
\end{figure}

Notice that the composition map
$h_{1_{l_{1}}}\circ f_{1}\circ h_{1_{k_{1}}}^{-1}$ is isotopic to the identity map 
$I_{(C,\tilde{K})}$ and this fact will be used to 
extend the map $f_{1}$ to a map  
$F_{1}:(V_{1}-\mbox{Int}( V_{2}))\rightarrow (V_{1}-\mbox{Int}( V_{2}))$ via the following Lemma.

\begin{lem}
Let $M$ be a compact 3-manifold with $\partial M\neq\emptyset$ (no necessarily connected). Let 
$g:\partial M\rightarrow \partial M$ be a homeomorphism which is
isotopic to the identity map $I_{\partial M}$. Then $g$ admits an extension $G:M\rightarrow M$. 
Furthermore,
suppose in addition that there exists a locally trivial fibration $\pi:M\to\mathbb S^1$ such that its 
restriction to $\partial M$ is also a locally trivial fibration and that $g$ leaves invariant the fibers
in the boundary. Then $g$ can be extended to a homeomorphism which preserves the fibers of $\pi$.
\end{lem}
{\it Proof.} Let $\psi:\partial M\times[0,1]\to M$ be a collaring of the boundary, i.e., 
$\psi$ is a homeomorphism such that $\psi(x,0)=x$. 
Let $N=\psi(\partial M\times[0,1])\subset M$.
Let $\{g_{t} \},\,\,t\in[0,1],$ be an isotopy
of $g$ to the identity, i.e., $g_o=g$, $g_1=I_{\partial M}$. Let 
$H:\partial M\times[0,1]\to \partial M\times[0,1]$
be given by the formula $H(x,t)=(g_t(x),t)$. Let $G_0:N\to N=    \psi\circ H\circ\psi^{-1}$.
Then we define $G:M\to M$ as $G(y)=G_0(y)$ if $y\in N$ and $G(y)=y$ if $y\notin N$. For the rest
we simply observe that the fibers of $\pi$ are surfaces that meet transversally the boundary
an therefore the collaring can be chosen in such a way that $\psi(\{(x,t)|t\in[0,1]\})$ is contained in
a fiber for each fixed $x\in\partial M$.
$\blacksquare$

We continue
inductively, so at the $k$-stage, we have a homeomorhism
$F_{k}:(V_{k}-\mbox{Int}( V_{k+1}))\rightarrow (V_{k}-\mbox{Int}( V_{k+1}))$ such that the solid cylinder
$(\partial V_{1_{i_{1}}2_{i_{2}}\ldots k_{i_{k}}}, \partial B_{1_{i_{1}}2_{i_{2}}\ldots k_{i_{k}}}\cap V_{k})$ 
is sent into 
$(\partial V_{1_{j_{1}}2_{j_{2}}\ldots k_{j_{k}}}, \partial B_{1_{j_{1}}2_{j_{2}}\ldots k_{j_{k}}}\cap V_{k})$, 
in such a way that the cylinder
$\partial V_{1_{i_{1}}2_{i_{2}}\ldots (k+1)_{i_{k+1}}}\cap V_{k+1}$ goes into the cylinder
$\partial V_{1_{j_{1}}2_{j_{2}}\ldots (k+1)_{j_{k+1}}}\cap V_{k+1}$ and 
($\pi_{1}^{-1}(p_{1_{i_{1}}2_{i_{2}}\ldots k_{i_{k}}})-\mbox{Int}(\pi_{2}^{-1}(p_{1_{i_{1}}2_{i_{2}}\ldots(k+1)_{i_{k+1}}}))$) 
is sent to
($\pi_{1}^{-1}(p_{1_{j_{1}}2_{j_{2}}\ldots k_{j_{k}}})-\mbox{Int}(\pi_{2}^{-1}(p_{1_{j_{1}}2_{j_{2}}\ldots(k+1)_{j_{k+1}}}))$).

This construction allows us to define a map 
$F:(\mathbb{S}^{3}-\Lambda)\rightarrow (\mathbb{S}^{3}-\Lambda)$
as $F(x)=F_{k}(x)$ if $x\in (V_{k}-\mbox{Int}(V_{k+1}))$. Notice that $F$ is a homeomorphism, since
each $F_{k}$ is a homeomorphism and $F_{k}(x)=F_{k+1}(x)$ for $x\in\partial V_{k+1}$ and for all $k$.
We extend $F$ to a map $\tilde{F}:\mathbb{S}^{3}\rightarrow\mathbb{S}^{3}$ in the following way. Let
$x\in\Lambda$. Then, there exists a sequence of cylinders  
$\{V_{1_{j_{1}},2_{j_{2}},\ldots,n_{j_{n}}}\}$, where
$V_{1_{j_{1}},2_{j_{2}},\ldots,n_{j_{n}}}\subset V_{n}$ such that 
$x=\cap_{n=1}^{\infty} V_{1_{j_{1}},2_{j_{2}},\ldots,n_{j_{n}}}$. We define 
$\tilde{F}(x)=\cap {F}(V_{1_{j_{1}},2_{j_{2}},\ldots,n_{j_{n}}})$. Notice that $\tilde{F}$ is well-defined 
and is continuous. In fact, since $F$ is a homeomorphism, we just need to prove
that $\tilde{F}$ is continuous in $\Lambda$. Given $x\in\Lambda$ and let $\{x_{n}\}$ be a sequence that 
converges to $x$.
We can assume, without loss of generality, that 
$x_{n}\in V_{1_{j_{1}},2_{j_{2}},\ldots,n_{j_{n}}}\subset V_{n}$, hence 
$\tilde{F}(x_{n})=F(x_{n})\in F(V_{j_{1},j_{2},\ldots,j_{n}})$, so 
$\lim_{n\rightarrow\infty}F(x_{n})=\tilde{F}(x)$. Therefore, $\tilde{F}$ is continuous.

\begin{theo}
The map $\tilde{F}:\mathbb{S}^{3}\rightarrow\mathbb{S}^{3}$ is a homeomorphism such that 
$\tilde{F}|_{\Lambda}=\Lambda$
and $\tilde{F}(p)=\tilde{F}(q)$.
\end{theo}
{\it Proof.} 
Since $\Lambda=\cap_{k=1}^{\infty}V_{k}$, and $\tilde{F}(V_{k})=V_{k}$, we have that 
$\tilde{F}(\Lambda)=\Lambda$.

Next, we will prove that $\tilde{F}$ is a bijection. Since $\tilde{F}|_{\mathbb{S}^{3}-\Lambda}$ 
is a bijection, it is
enough to prove that $\tilde{F}|_{\Lambda}$ is it.
Let $a$, $b\in\Lambda$. Then $a=\cap V_{1_{l_{1}},2_{l_{2}},\ldots,n_{l_{n}}}$ and 
$b=\cap V_{1_{r_{1}},2_{r_{2}},\ldots,n_{r_{n}}}$, where 
$V_{1_{l_{1}},2_{l_{2}},\ldots,n_{l_{n}}}$, $V_{1_{r_{1}},2_{r_{2}},\ldots,n_{r_{n}}}\subset V_{n}$. If 
$\tilde{F}(a)=\tilde{F}(b)$, this implies that 
$\tilde{F}(V_{1_{l_{1}},2_{l_{2}},\ldots,n_{l_{n}}})\cap 
\tilde{F}(V_{1_{r_{1}},2_{r_{2}},\ldots,n_{r_{n}}})\neq\emptyset$.
Since, 
$\tilde{F}|_{\mathbb{S}^{3}-\Lambda}$ is a homeomorphism, we have that 
$V_{1_{l_{1}},2_{l_{2}},\ldots,n_{l_{n}}}\cap V_{1_{r_{1}},2_{r_{2}},\ldots,n_{r_{n}}}\neq\emptyset$, 
but this is a contradiction. Hence $a=b$. 
For each $x=\cap V_{1_{j_{1}},2_{j_{2}},\ldots,n_{j_{n}}}\in\Lambda$,  let 
$x'=\cap \tilde{F}^{-1}(V_{1_{j_{1}},2_{j_{2}},\ldots,n_{j_{n}}})$. By the above, 
$\tilde{F}(x')=x$. It follows that, $\tilde{F}|_{\Lambda}:\Lambda\rightarrow\Lambda$ is a bijection, 
hence $\tilde{F}$
is a bijection.
Therefore $\tilde{F}:\mathbb{S}^{3}\rightarrow \mathbb{S}^{3}$ is a continuous bijection, hence 
$\tilde{F}$ is a homeomorphism. $\blacksquare$

\begin{coro}
Dynamically-defined wild knots are homogeneous.
\end{coro} 
{\it Proof.}
Let $p$, $q\in \Lambda$. Then $p=\cap V_{1_{l_{1}},2_{l_{2}},\ldots,n_{l_{n}}}$ and 
$q=\cap V_{1_{r_{1}},2_{r_{2}},\ldots,n_{r_{n}}}$, where 
$V_{1_{l_{1}},2_{l_{2}},\ldots,n_{l_{n}}}$, $V_{1_{r_{1}},2_{r_{2}},\ldots,n_{r_{n}}}\in V_{n}$. 
By the above, we can construct a homeomorphism
$\tilde{F}:\mathbb{S}^{3}\rightarrow\mathbb{S}^{3}$ such that $\tilde{F}|_{\Lambda}=\Lambda$
and $\tilde{F}(p)=\tilde{F}(q)$. Therefore, $\Lambda$ is homogeneous. $\blacksquare$

\begin{rem}
The same method applies to prove the following theorem.
\begin{theo}
Let $T_{k}\subset\mathbb{S}^{3}$ be a nested decreasing sequence of smooth solid tori ($k\in\mathbb{N}$),
i.e. $T_{k+1}\subset \mbox{Int}(T_{k})$. Suppose that $\cap^{\infty}_{k=1}T_{k}:=K$ is a knot
(wild or not). Then $K$ is homogeneous.
\end{theo}
\end{rem}

\begin{rem}
Let $\gamma\subset\Lambda$ be an 
orbit under the Kleinian group $\Gamma$. Notice that the action of $\Gamma$ is minimal, i.e. the
orbits are dense. Let $x$, $y\in\gamma$.
Then, using the action of $\Gamma$, we can construct a homeomorphism 
$H:\mathbb{S}^{3}\rightarrow\mathbb{S}^{3}$ such that $H|_{\Lambda}=\Lambda$
and $H(p)=H(q)$. However, the above theorem holds for any couple of points $p$, $q\in\Lambda$. 
\end{rem}
\section{Dynamically-defined fibered wild knots}

We recall that a knot or link $L$ in $\mathbb{S}^{3}$ is {\it fibered} if there exists a 
locally trivial fibration $f:(\mathbb{S}^{3}-L)\rightarrow \mathbb{S}^{1}$. We 
require that $f$ be well-behaved near $L$, that is, each component $L_{i}$ is 
to have a neighbourhood framed as $\mathbb{D}^{2}\times\mathbb{S}^{1}$, with 
$L_{i}\cong \{0\}\times\mathbb{S}^{1}$, in such a way that the
restriction of $f$ to $(\mathbb{D}^{2}-\{0\})\times\mathbb{S}^{1}$ is the map 
into $\mathbb{S}^{1}$ given by $(x,y)\rightarrow \frac{y}{|y|}$. It follows that each
$f^{-1}(x)\cup L$, $x\in\mathbb{S}^{1}$, is a 2-manifold
 with boundary $L$: in fact a Seifert surface for $L$ 
(see \cite{rolfsen}, page 323).

\begin{exs}
The right-handed trefoil knot and the figure-eight knot are fibered knots with fiber
the punctured torus.
\end{exs}

For dynamically-defined wild knots we have the following result.
\begin{theo}
Let $\Sigma_{1}$, $\Sigma_{2},\ldots ,\Sigma_{n}$ be round 2-spheres
in $\mathbb{S}^{3}$ which form a necklace whose template is a non-trivial
tame fibered knot $K$. Let $\Gamma$ be the group generated by reflections
$I_{j}$ on $\Sigma_{j}$ ($j=1,2,\ldots,n$) and let $\widetilde{\Gamma}$
be the orientation-preserving index two subgroup of $\Gamma$. 
Let $\Lambda(\Gamma)=\Lambda(\widetilde{\Gamma})$ be the 
corresponding limit set. Then:
\begin{enumerate}
\item There exists a locally trivial fibration $\psi
:(\mathbb{S}^{3}-\Lambda(\Gamma))\rightarrow\mathbb{S}^{1}$, where the
 fiber $\Sigma_{\theta}=\psi^{-1}(\theta)$ is an orientable
 infinite-genus surface with one end.
\item  $\overline{\Sigma_{\theta}}-\Sigma_{\theta}=\Lambda(\Gamma)$, where
$\overline{\Sigma_{\theta}}$ is the closure of $\Sigma_{\theta}$ in $\mathbb{S}^{3}$.
\end{enumerate}
\end{theo}

Next, we will briefly describe the fiber $\Sigma_{\theta}$. For a proof 
of the above theorem, see \cite{gaby}.  

Let $T_{1}$ be a pearl-necklace subordinate to the fibered tame knot $K_{1}$ with fiber $S_{1}$. 
Let $\widetilde{P}:(\mathbb{S}^{3}-K_{1})\rightarrow \mathbb{S}^{1}$ be the
given fibration. Observe that 
$\widetilde{P}\mid_{\mathbb{S}^{3}-|T_{1}|}\equiv P$ is a fibration and, after
modifying $\widetilde{P}$ by isotopy if necessary, we can consider  that
the fiber $S$  cuts each pearl $\Sigma_{i}\in T_{1}$ in arcs
$a_{i}$, whose end-points are $\Sigma_{i-1}\cap\Sigma_{i}$ and $\Sigma_{i}\cap\Sigma_{i+1}$. These two
points belong to the limit set (see Figure 9).

\begin{figure}[tbh]
\centerline{\epsfxsize=1.2in \epsfbox{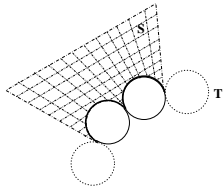}}
\caption{\sl The fiber intersects each pearl in arcs.}
\label{F10}
\end{figure}

The fiber $\widetilde{P}^{-1}(\theta)=P^{-1}(\theta)$ is a Seifert surface $S_{1}^{*}$ of $K_{1}$, 
for each $\theta\in\mathbb{S}^{1}$. We suppose $S_{1}^{*}$ is oriented.

The reflection  $I_{j}$ maps both a copy of $T_{1}-\Sigma_{j}$
(called $T_{1}^{j}$) and a copy of $S_{1}^{*}$ (called $S_{2}^{*j}$) into
the ball $B_{j}$, for $j=1,2,\ldots,n$. Observe that both $T_{1}^{j}$ and $S_{2}^{*j}$ 
have opposite orientation 
and that $S_{1}^{*}$ and $S_{2}^{*j}$ are joined by
the arc $a_{j}$ (see Figure 10) which, in both surfaces, has the same orientation. 

\begin{figure}[tbh]
\centerline{\epsfxsize=.8in \epsfbox{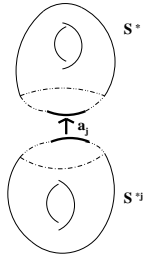}}
\caption{\sl Sum of two surfaces $S_{1}^{*}$ and $S_{2}^{*j}$ along arc $a_{j}$.}
\label{F11}
\end{figure}

The necklaces $T_{1}^{j}$ and $T_{1}$ are joined by the points of tangency 
of the pearl $\Sigma_{j}$ and the 
orientation of these two points is preserved by the reflection $I_{j}$. Thus, we have obtained a new 
pearl-necklace isotopic to the 
connected sum $T_{1}\#T_{1}^{j}$, whose complement also fibers over the circle with 
fiber the sum of $S_{1}^{*}$ with $S_{2}^{*j}$ along arc $a_{j}$, namely
the fiber is $S_{1}^{*}\#_{a_{j}}S_{2}^{*j}$.

Now do this for each $j=1,\ldots,n$. At the end of the first stage, we
have a new pearl-necklace $T_{2}$ whose template is the knot $K_{2}$
(see section 2). Its complement fibers 
over the circle with fiber the Seifert surface $S^{*}_{2}$, which is in turn 
homeomorphic to the sum of $n+1$ copies of $S_{1}^{*}$ along the respective arcs.

Continuing this process, we have from the second step onwards, that $n-1$ copies of $S_{1}^{*}$ are added 
along arcs to 
each surface $S^{*i}_{k}$, (the surface corresponding 
to the $k^{th}$ stage). Notice that 
in each step,
the points of tangency are removed since they belong to the limit set,
and the length of the arcs $a_{j}$ tends to zero.

From the remarks above, we have that $\Sigma_{\theta}$ is homeomorphic to an orientable
infinite-genus surface. In fact, it is the sum along arcs of an infinite
number of copies of $S^{*}$. 

Using the second part of lemma 3.1 one has the following:

\begin{theo}
Let $K_{1}$ be a tame fibered knot with fiber $S_{1}$. Let $\Lambda$ be the wild knot 
obtained from $K_{1}$ through
a reflecting process. Then, given $p$, $q\in\Lambda$, 
there exists a homeomorphism $\tilde{F}:\mathbb{S}^{3}\rightarrow\mathbb{S}^{3}$ such that 
$\tilde{F}|_{\Lambda}=\Lambda$, $\tilde{F}(p)=\tilde{F}(q)$ and  
preserves the fibers.
\end{theo}
{\it Proof.}
Let $P:(\mathbb{S}^{3}-K)\rightarrow \mathbb{S}^{1}$ be the given fibration. We can assume that
the fiber cuts each pearl $\Sigma_{i}$ of the pearl-necklace as in Figure 9. 
Let $p$, $q\in\Lambda$. Then, there exist two sequences of solid cylinders  
$\{V_{1_{i_{1}},2_{i_{2}},\ldots,n_{i_{n}}}\}$,  and  
$\{V_{1_{j_{1}},2_{j_{2}},\ldots,n_{j_{n}}}\}$ where $V_{1_{i_{1}},2_{i_{2}},\ldots,n_{i_{n}}}$ 
and $V_{1_{j_{1}},2_{j_{2}},\ldots,n_{j_{n}}}\in V_{n}$, such that 
$p=\cap_{n=1}^{\infty} V_{1_{i_{1}},2_{i_{2}},\ldots,n_{i_{n}}}$ and
  $q=\cap_{n=1}^{\infty} V_{1_{j_{1}},2_{j_{2}},\ldots,n_{j_{n}}}$. Let $V_{k}$ be as in section 3. 
Consider the map 
$f_{k}:(\partial V_{k}\cap\partial V_{k+1})\rightarrow (\partial V_{k}\cap \partial V_{k+1})$.

By the previous description, we know that $V_{1}-\mbox{Int}(V_{2})$ fibers over the circle and that $f_{1}$
leaves invariant the fibers on the boundary. Then by Lemma 3.1, $f_{1}$ admits and extension 
$F_{1}:(V_{1}-\mbox{Int}(V_{2}))\rightarrow (V_{1}-\mbox{Int}(V_{2}))$ which preserves the fibers. 
We continue inductively, so at the end of the $k$-stage, we have that the homeomorphism $f_{k}$ preserves 
the fibers, hence it can be extended to a map 
$F_{k}:(V_{k}-\mbox{Int}(V_{k+1}))\rightarrow (V_{k}-\mbox{Int}(V_{k+1}))$ which 
also preserves the fibers. 

Let $\tilde{F}:\mathbb{S}^{3}\rightarrow \mathbb{S}^{3}$ be
as in section 3. Then, by Theorem 3.2 we have that 
$\tilde{F}|_{\Lambda}=\Lambda$, $\tilde{F}(p)=\tilde{F}(q)$. The map 
$F=\tilde{F}|_{\mathbb{S}^{3}-\Lambda}$ defined by $F(x)=F_{k}(x)$ if 
$x\in (V_{k}-\mbox{Int}(V_{k+1}))$ is a homeomorphism
(see section 3) which preserves fibers. Therefore, the result follows. $\blacksquare$

G. Hinojosa. Universidad Aut\'oma del Estado de Morelos. Av. Universidad 1001, Col. Chamilpa.
Cuernavaca, Morelos, M\'exico, 62210. 

\hspace{-.6cm}{\it E-mail address:}gabriela@servm.fc.uae.mx
\vskip .5cm
A. Verjovsky. Instituto de Matem\'aticas, Universidad Nacional Au\-t\'o\-no\-ma de M\'exico, Unidad
Cuernavaca. Av. Universidad s/n, Col. Lomas de Chamilpa. Cuernavaca, Morelos M\'exico, 62210.

\hspace{-.6cm}{\it E-mail address:} alberto@matcuer.unam.mx 

\end{document}